\documentclass[twoside]{ima} 

\usepackage{clrscode}

\usepackage{graphicx,subfigure,psfrag,amsmath,amssymb}
\usepackage[dvips]{color}
\usepackage{lscape,longtable}

\usepackage{rotating}

\usepackage[dvips]{geometry}

\newcommand{\st}{\mbox{s.t.}}

\newcommand{\SDP}{$PSD$}
\newcommand{\bSDP}{{\bf PSD}}
\newcommand{\bSDPLP}{{\bf PSDLP}}
\newcommand{\SDPCUTS}{PSDCUT}
\newcommand{\RLT}{$RLT$}

\begin{document}

\markboth{ANDREA QUALIZZA ET AL.\mbox{}}{LP RELAXATIONS
  OF QUADRATIC PROGRAMS}

\title{LINEAR PROGRAMMING RELAXATIONS OF QUADRATICALLY CONSTRAINED
  QUADRATIC PROGRAMS}

\author{
ANDREA QUALIZZA${}^{1}$ \and 
PIETRO BELOTTI${}^{2}$ \and
FRAN\c{C}OIS MARGOT${}^{1,3}$
}

\maketitle   

\footnotetext[1] {
Tepper School of Business, Carnegie Mellon University, Pittsburgh PA.}

\footnotetext[2] {
Dept.\mbox{} of Industrial and Systems Engineering, 
Lehigh University, Bethlehem PA. }

\footnotetext[3] {
Corresponding author (email: {\tt fmargot@andrew.cmu.edu}).
Supported by NSF grant NSF-0750826.
}%

\begin{abstract}
We investigate the use of linear programming tools for solving semidefinite
programming relaxations of quadratically constrained quadratic problems.
Classes of valid linear inequalities are presented, including
sparse \SDP{} cuts, and principal minors \SDP{} cuts.
Computational results based on instances from the literature are presented.
\end{abstract}             
 

{\AMSMOS 
90C57
\endAMSMOS}

\section{Introduction}
\label{introduction}

Many combinatorial problems have Linear Pr\-ogr\-amm\-ing ({\it LP}) relaxations
that are commonly used for their solution through branch-and-cut algorithms.
Some of them also have stronger relaxations involving positive semidefinite
(\SDP{}) constraints. In general, stronger relaxations should
be preferred when solving a problem, thus using these \SDP{} relaxations
is tempting. However, they come with the drawback of requiring a 
Semidefinite Programming (SDP)
solver, creating practical difficulties for an efficient implementation
within a branch-and-cut algorithm. Indeed, a major weakness of current 
SDP solvers 
compared to LP solvers is their lack of efficient warm starting mechanisms.
Another weakness is solving problems involving a mix of \SDP{} constraints
and a large number of linear inequalities, as these linear inequalities
put a heavy toll on the linear algebra steps required during the solution 
process.

In this paper, we investigate LP relaxations of \SDP{} constraints with the
aim of capturing most of the strength of the \SDP{} relaxation, while
still being able to use an LP solver. The LP relaxation we obtain is an 
outer-approximation of the \SDP{} cone, with the typical convergence
difficulties when aiming to solve problems to optimality. We thus
do not cast this work as an efficient way to solve \SDP{} problems,
but we aim at finding practical ways to approximate \SDP{} constraints
with linear ones. 

We restrict our experiments to Quadratically Constrained Quadratic Program
({\it QCQP}). A QCQP problem with variables 
$x \in \mathbb{R}^n$ and $y \in \mathbb{R}^m$ is a problem of the form
\begin{displaymath}
\tag*{({\bf QCQP})}
\begin{array}{lllll}
\max 	& x^T Q_0 x + a_0^T x + b_0^T y \\
\st \\
	& x^T Q_k x + a_k^T x + b_k^T y & \leq c_k 
& \mbox{for \ } k = 1, 2, \ldots, p \\
	& l_{x_i} \leq x_i \leq u_{x_i}	  &	     
& \mbox{for \ } i = 1,2,\ldots,n\\
	& l_{y_j} \leq y_j \leq u_{y_j}	  &	     
& \mbox{for \ } j = 1,2,\ldots,m
\end{array}
\end{displaymath}

\noindent
where, for $k = 0, 1, 2, \ldots, p$,
$Q_k$ is a rational symmetric $n \times n$-matrix,
$a_k$ is a rational $n$-vector,
$b_k$ is a rational $m$-vector, and $c_k \in \mathbb{Q}$.
Moreover, the lower and upper bounds $l_{x_i}, u_{x_i}$ for $i = 1, \ldots, n$,
and $l_{y_j}, u_{y_j}$ for $j = 1, \ldots, m$ are all finite.
If $Q_0$ is negative semidefinite and $Q_k$ is positive semidefinite for each 
$k = 1, 2, \ldots, p$, problem {\bf QCQP} is convex and 
thus easy to solve. Otherwise, the problem is NP-hard \cite{BurerLetchford}. 

An alternative {\it lifted} formulation for {\bf QCQP} is obtained by 
replacing each  quadratic term $x_i x_j$ with a new variable $X_{ij}$. 
Let $X = x x^T$ be the matrix with entry $X_{ij}$ 
corresponding to the quadratic term $x_i x_j$. For square matrices $A$ and $B$
of the same dimension, let $A \bullet B$ denote the {\it Frobenius inner
product} of $A$ and $B$, i.e., the trace of $A^T B$.
Problem {\bf QCQP} is then equivalent to 

\begin{displaymath}
\tag*{({\bf LIFT})}
\begin{array}{lllll}
\max 	& Q_0 \bullet X + a_0^T x + b_0^T y \\
\st \\
	& Q_k \bullet X + a_k^T x + b_k^T y & \leq c_k & 
\mbox{for \ } k = 1, 2, \ldots, p \\
	& l_{x_i} \leq x_i \leq u_{x_i}	  &	     & 
\mbox{for \ } i = 1,2,\ldots,n\\
	& l_{y_j} \leq y_j \leq u_{y_j}	  &	     & 
\mbox{for \ } j = 1,2,\ldots,m\\
	& X = x x^T.
\end{array}
\end{displaymath}

The difficulty in solving problem {\bf LIFT} lies in the non-convex constraint
$X = xx^T$. A relaxation, dubbed \SDP{}, that is possible to solve relatively 
efficiently is obtained by relaxing 
this constraint to the requirement that $X - x x^T$ be positive 
semidefinite, i.e., $X - x x^T \succeq 0$. An alternative relaxation 
of {\bf QCQP}, dubbed \RLT{}, is obtained by the Reformulation Linearization 
Technique \cite{SheraliAdams}, using products of pairs of 
original constraints and bounds and replacing nonlinear terms with new 
variables.

Anstreicher \cite{Anstreicher} compares the \SDP{} and \RLT{} relaxations 
on a set of quadratic problems with box constraints, i.e., {\bf QCQP}
problems with $p=0$ and with all the variables bounded between 
0 and 1. He shows that the \SDP{} relaxations of these instances are fairly 
good and that combining the \SDP{} and \RLT{} relaxations yields significantly
tighter relaxations than either of the \SDP{} or \RLT{} relaxations. 
The drawback of combining the two relaxations is that current 
SDP solvers have difficulties to handle the large 
number of linear constraints of the \RLT{}. 

Our aim is to solve relaxations of {\bf QCQP} using exclusively linear 
programming tools. The \RLT{} is readily applicable for our purposes, 
while the \SDP{} technique requires a cutting plane approach as described 
in Section~\ref{relaxations}.

In Section~\ref{framework} we consider several families of valid cuts. 
The focus is essentially on capturing the strength of the positive
semidefinite condition using standard cuts \cite{Fraticelli}, and some 
sparse versions of these. 

We analyze empirically the strength of the considered cuts on instances 
taken from GL\-OB\-ALLib \cite{GLOBALLib} and quadratic programs with box 
constraints described in more details in the next section. 
Implementation and computational results are presented in 
Section~\ref{computationalresults}.
Finally, Section~\ref{conclusions} summarizes the results and gives 
possible directions for future research.

\section{Relaxations of QCQP problems}
\label{relaxations}
A typical approach to get bounds on the optimal value of a {\bf QCQP} is 
to solve a convex relaxation. Since our aim is to work with linear relaxations, 
the first step is to linearize {\bf LIFT} by relaxing the last constraint
to $X = X^T$. We thus get the Extended formulation
\begin{displaymath}
\tag*{({\bf EXT})}
\begin{array}{lllll}
\max 	& Q_0 \bullet X + a_0^T x + b_0^T y \\
\st \\
	& Q_k \bullet X + a_k^T x + b_k^T y & \leq c_k & 
\mbox{for \ } k = 1, 2, \ldots, p \\
	& l_{x_i} \leq x_i \leq u_{x_i}	  &	     & 
\mbox{for \ } i = 1,2,\ldots,n\\
	& l_{y_j} \leq y_j \leq u_{y_j}	  &	     & \mbox{for \ } 
j = 1,2,\ldots,m\\
	& X = X^T.
\end{array}
\end{displaymath}
{\bf EXT} is a Linear Program with $n(n+3)/2 + m$ variables and the same 
number of constraints as {\bf QCQP}. Note that the optimal value of {\bf EXT}
is usually a weak upper bound for {\bf QCQP}, as no constraint links
the values of the $x$ and $X$ variables. Two main approaches for doing that
have been proposed and are based on relaxations of the last constraint of
{\bf LIFT}, namely
\begin{equation}
X - x x^T = 0 \label{X-xxTeq0} .
\end{equation}

They are known as the Positive Semidefinite (\SDP{}) relaxation
and the Reformulation Linearization Technique (\RLT{}) relaxation.

\subsection{PSD Relaxation}
\label{sdprelaxation}
As $X - x x^T = 0$ implies  $X - x x^T \succcurlyeq 0$, using this last
constraint yields a convex relaxation of {\bf QCQP}. 
This is the approach used in 
\cite{Fraticelli,KrishnanMitchell1,KrishnanMitchell2,VandenbergheBoyd}, 
among others.

Moreover, using Schur's complement
\begin{displaymath}
X - x x^T \succcurlyeq 0
\hspace{1cm}
\Leftrightarrow 
\hspace{1cm}
\begin{pmatrix}
1 & x^T\\
x & X
\end{pmatrix} \succcurlyeq 0,
\end{displaymath}
and defining
\begin{displaymath}
\tilde{Q}_k = 
\begin{pmatrix}
-c_k & a_k^T/2\\
a_k/2 & Q_k
\end{pmatrix} ,
\hspace{0.5in}  
\tilde{X} = 
\begin{pmatrix}
1 & x^T\\
x & X
\end{pmatrix},
\end{displaymath}
we can write the \SDP{} relaxation of {\bf QCQP} in the compact form
\begin{displaymath}
\tag*{(\bSDP)}
\begin{array}{llll}
\max 	& \tilde{Q}_0 \bullet \tilde{X} + b_0^T y \\
\st \\
	& \tilde{Q} \bullet \tilde{X} + b_k^T y  \leq 0	&& 
k = 1, 2, \ldots, p \\
	& l_{x_i} \leq x_i \leq u_{x_i}	  && i = 1,2,\ldots,n\\
	& l_{y_j} \leq y_j \leq u_{y_j}	  && j = 1,2,\ldots,m\\
	& \tilde{X}  \succcurlyeq 0.
\end{array}
\end{displaymath}
This is a positive semidefinite problem with linear constraints.
It can thus be solved in polynomial time using interior point algorithms.
\bSDP{} is tighter than usual linear relaxations for problems such as
the Maximum Cut, Stable Set, and Quadratic Assignment problems 
\cite{Wolkowicz}. All these problems can be formulated as {\bf QCQP}s.

\subsection{RLT Relaxation}
\label{rltrelaxation}
The Reformulation Linearization Technique \cite{SheraliAdams} can be used
to produce a relaxation of {\bf QCQP}. It adds linear inequalities to
{\bf EXT}. These inequalities are derived from the variable bounds and 
constraints of the original problem as follows: multiply together 
two original constraints or bounds and replace each product term 
$x_i x_j$ with the variable $X_{ij}$.
For instance, let $x_i, x_j,\ i,j\in \{1,2,\ldots,n\}$ be two variables from 
{\bf QCQP}. By taking into account only the four original bounds 
$x_i - l_{x_i} \geq 0$, $x_i - u_{x_i} \leq 0$, $x_j - l_{x_j} \geq 0$, 
$x_j - u_{x_j} \leq 0$, we get the $RLT$ inequalities
\begin{equation}
\label{rltcutsfrombounds}
\begin{array}{llllllll}
X_{ij}  -   l_{x_i} x_j  -  l_{x_j} x_i  & \geq  -  l_{x_i} l_{x_j}, \\
X_{ij}  -   u_{x_i} x_j  -  u_{x_j} x_i  & \geq  -  u_{x_i} u_{x_j}, \\
X_{ij}  -   l_{x_i} x_j  -  u_{x_j} x_i  & \leq  -  l_{x_i} u_{x_j}, \\
X_{ij}  -   u_{x_i} x_j  -  l_{x_j} x_i  & \leq  -  u_{x_i} l_{x_j} \ .
\end{array}
\end{equation}

Anstreicher \cite{Anstreicher} observes that, for Quadratic Programs 
with box constraints, the \SDP{} and \RLT{} constraints together yield 
much better bounds than those obtained from the \bSDP{} or {\bf RLT} 
relaxations. In this work, we want to capture the strength of both 
techniques and generate a Linear Programming relaxation of {\bf QCQP}.

Notice that the four inequalities above, introduced by McCormick
\cite{mccormick}, constitute the convex envelope of the set
$\{(x_i,x_j,X_{ij}) \in \mathbb R^3: 
l_{x_i} \le x_i \le u_{x_i}, 
l_{x_j} \le x_j \le u_{x_j}, 
X_{ij} = x_i x_j\}$ as proven by Al-Khayyal and
Falk \cite{al-khayyal}, i.e., they are the tightest relaxation for the
single term $X_{ij}$. 

\section{Our Framework}
\label{framework}
While the \RLT{} constraints are linear in the variables in the {\bf EXT} 
formulation and therefore can be added directly to {\bf EXT}, this is not 
the case for the \SDP{} constraint. We use a linear outer-approximation 
of the \bSDP{} relaxation and a cutting plane framework, adding 
a linear inequality separating the current solution from the \SDP{} cone.

The initial relaxation we use and the various cuts generated
by our separation procedure are described
in more details in the next sections.

\subsection{Initial Relaxation}
\label{initialrelaxation}
Our initial relaxation is the {\bf EXT} formulation together with the $O(n^2)$ 
\RLT{} constraints derived from the bounds on the variables $x_i,\ i=1,2,\ldots,n$. 
We did not include the \RLT{} constraints derived from the problem 
constraints due to their large number and the fact that we want to avoid
the introduction of extra variables for the multivariate 
terms that occur when quadratic constraints are multiplied together.

The bounds $[L_{ij},U_{ij}]$ for the extended variables $X_{ij}$ are 
computed as follows:
\begin{displaymath}
\begin{array}{ll}
L_{ij}  =  \min  \{ l_{x_i} l_{x_j};\,  l_{x_i} u_{x_j};\,  u_{x_i} l_{x_j};\,  u_{x_i} u_{x_j}\},\!\!&\forall  i=1,2,\ldots,n;\ j=i,\ldots,n\\
U_{ij}  =  \max  \{ l_{x_i} l_{x_j};\,  l_{x_i} u_{x_j};\,  u_{x_i} l_{x_j};\,  u_{x_i} u_{x_j}\},\!\!&\forall  i=1,2,\ldots,n;\ j=i,\ldots,n.
\end{array}
\end{displaymath}
In addition, equality (\ref{X-xxTeq0}) implies $X_{ii} \geq x_i^2$.
We therefore also make sure that $L_{ii} \ge 0$. 
In the remainder of the paper, this initial relaxation is identified as
{\bf EXT+RLT}.

\subsection{PSD Cuts}
\label{sdpcuts}

We use the equivalence that a matrix is positive semidefinite
if and only if 
\begin{equation}
\label{SDPlinearcons}
v^T \tilde{X} v \geq 0 \quad \mbox{for \ all \ } v\in \mathbb{R}^{n+1} \ .
\end{equation}
We can reformulate \bSDP{} as the semi-infinite Linear Program
\begin{displaymath}
\tag*{(\bSDPLP{})}
\begin{array}{llll}
\max 	& \tilde{Q}_0 \bullet \tilde{X} + b_0^T y \\
\st \\
	& \tilde{Q} \bullet \tilde{X} + b_k^T y & \leq c_k & 
\mbox{for \ } k = 1, 2, \ldots, p \\
	& l_{x_i} \leq x_i \leq u_{x_i}	  &	     & 
\mbox{for \ }  i = 1,2,\ldots,n\\
	& l_{y_j} \leq y_j \leq u_{y_j}	  &	     & 
\mbox{for \ } j = 1,2,\ldots,m\\
	& v^T \tilde{X} v     & \geq 0	     & 
\mbox{for \ all \ }  v\in \mathbb{R}^{n+1}.
\end{array}
\end{displaymath}

A practical way to use \bSDPLP{} is to adopt a cutting plane approach to 
separate constraints (\ref{SDPlinearcons}) as done in
\cite{Fraticelli}.

Let $\tilde{X}^*$ be an arbitrary point in the space of the $\tilde{X}$ variables. 
The spectral decomposition of $\tilde{X}^*$ is used to decide if  $\tilde{X}^*$
is in the \SDP{} cone or not. Let the eigenvalues and corresponding 
orthonormal eigenvectors of $\tilde{X}^*$ be $\lambda_k$ and $v_k$ for 
$k = 1,2,\ldots,n$, and assume without loss of generality that 
$\lambda_1 \le \lambda_2 \le \ldots \le \lambda_n$ and let $t \in \{0, \ldots, n\}$ 
such that $\lambda_t < 0 \le \lambda_{t+1}$.
If $t = 0$, then all the eigenvalues are non negative and
$\tilde{X}^*$ is positive semidefinite. Otherwise, 
$v_k^T \tilde{X}^* v_k = v_k^T \lambda_k v_k = \lambda_k < 0$ for $k = 1, \ldots, t$. 
Hence, the valid cut
\begin{equation}
\label{CUTS:sdpcuts}
v_k^T \tilde{X} v_k \geq 0
\end{equation}
is violated by 
$\tilde{X}^*$. Cuts of the form (\ref{CUTS:sdpcuts}) are called \SDPCUTS{}s
in the remainder of the paper.

The above procedure has two major weaknesses: First, only one cut
is obtained from eigenvector $v_k$ for $k = 1, \ldots, t$, while computing
the spectral decomposition requires a non trivial investment in cpu time, 
and second, the cuts are usually very dense, i.e.\mbox{} almost all 
entries in $v v^T$ are nonzero. Dense cuts are frowned upon when used in
a cutting plane approach, as they might slow down considerably the 
reoptimization of the linear relaxation. 

To address these weaknesses, we describe in the next section
a heuristic to generate several
sparser cuts from each of the vectors  $v_k$ for $k = 1, \ldots, t$.

\subsection{Sparsification of PSD cuts}
\label{sparsify}

A simple idea to get sparse cuts is to start with vector $w = v_k$,
for $k = 1, \ldots, t$, and iteratively set to zero some component of $w$,
provided that $w^T \tilde{X}^* w$ remains sufficiently negative. 
If the entries are considered
in random order, several cuts can be obtained from a single eigenvector $v_k$.
For example, consider the {\em Sparsify} procedure in Figure \ref{FIG:sparsify},
 taking as parameters
an initial vector $v$, a matrix $\tilde{X}$, and two numbers between 0 and 1,
$pct_{NZ}$ and $pct_{VIOL}$, that control the maximum percentage  of
nonzero entries in the final vector and the minimum violation requested for
the corresponding cut, respectively. In the procedure, parameter $length[v]$
identifies the size of vector $v$.

\begin{figure}
\begin{codebox}
\Procname{$\proc{\textit{Sparsify}}(v,\tilde{X}, pct_{NZ}, pct_{VIOL})$}
\li $\id{min_{VIOL}} \gets - v^T \tilde{X} v \cdot pct_{VIOL}$
\li $\id{max_{NZ}} \gets \lfloor length[v] \cdot pct_{NZ} \rfloor$
\li $w \gets v$
\li $\id{perm} \gets \textrm{random permutation of 1 }\To \id{length}[w]$
\li \For $j \gets 1$ \To $\id{length}[w]$ \label{divers1} 
\li \Do
\li      $z \gets w$, $z[perm[j]] \gets 0$ \label{assignz}
\li	 \If $- z^T \tilde{X} z > min_{VIOL}$ \label{violz}
\li	 \Then	$w \gets z$ \label{w_i=0}
	 \End
    \End
\li \If $\mbox{number of non-zeroes in \ } w < max_{NZ}$
\li	\Then
		$\mbox{output \ } w$
	\End
\end{codebox}
\caption{Sparsification procedure for PSD cuts}
\label{FIG:sparsify}
\end{figure}

It is possible to implement this procedure to run in 
$O(n^2)$ if $length[v] = n+1$: Compute and update a vector $m$ such that
\begin{eqnarray*}
m_j= \sum_{i=1}^{n+1} w_j w_i \tilde{X}_{ij} &\mbox{ for \ }
j = 1, \ldots, n+1 \ .
\end{eqnarray*}
\noindent
Its initial computation takes $O(n^2)$ and
its update, after a single entry of $w$ is set to 0, takes $O(n)$. The vector
$m$ can be used to compute the left hand side of the test in step \ref{violz} 
in constant time given the value
of the violation $d$ for the inequality generated by the current vector $w$:
Setting the entry $\ell = perm[j]$ of $w$ to zero reduces the violation
by $2 m_\ell - w_\ell^2 \tilde{X}_{\ell\ell}$ 
and thus the violation of the resulting vector is $(d - 2 m_{\ell} + w_{\ell}^2
 \tilde{X}_{\ell\ell})$.

A slight modification of the procedure is used to obtain several cuts
from the same eigenvector: Change the loop condition in step \ref{divers1}
to consider the entries in $perm$ in cyclical order, from all possible starting
points $s$ in $\{1,2 \ldots, length[w]\}$, with the additional condition that 
entry $s-1$ is not set to 0 when starting from $s$ to guarantee that we do not
generate always the same cut.
From our experiments, this simple idea produces collections of sparse and 
well-diversified cuts. 
This is referred to as SPARSE1 in the remainder of the paper.

We also consider the following variant of the procedure given in 
Figure~\ref{FIG:sparsify}. 
Given a vector $w$, let $\tilde{X}_{[w]}$ be the principal minor of $\tilde{X}$
induced by the indices of the nonzero entries in $w$. Replace step 
\ref{assignz} with
\begin{quote}
  \ref{assignz}. $z \gets \bar w$ where $\bar w$ is an eigenvector
  corresponding to the most negative eigenvalue of a spectral
  decomposition of $\tilde{X}_{[w]}$, $z[perm[j]] \gets 0$.
\end{quote}
This is referred to as SPARSE2 in the remainder, and 
we call the cuts generated by SPARSE1 or SPARSE2 described above {\it Sparse 
\SDP{} cuts}.

Once sparse \SDP{} cuts are generated, for each vector $w$ generated,
we can also add
all \SDP{} cuts given by the eigenvectors corresponding to negative 
eigenvalues of a spectral decomposition of $\tilde{X}_{[w]}$. 
These cuts are valid and sparse. They are called 
{\it Minor} \SDP{} {\it cuts} and denoted by MINOR in the following.

An experiment to determine good values for the parameters $pct_{NZ}$ and 
$pct_{VIOL}$ was performed on the 38 GLOBALLIB instances and 51 BoxQP 
instances described
in Section~\ref{instances}. It is run by selecting two sets of three
values in $[0, 1]$, $\{V_{LOW},V_{MID},V_{UP}\}$ for $pct_{VIOL}$ and 
$\{N_{LOW},N_{MID},N_{UP}\}$ for $pct_{NZ}$.
The nine possible combinations of these parameter values
are used and the best of the nine $(V_{best},N_{best})$ is selected. We then 
center and reduce the possible ranges around $V_{best}$ and $N_{best}$, 
respectively, and repeat the operation. The procedure is stopped when the 
best candidate parameters are $(V_{MID},N_{MID})$ and the size of the 
ranges satisfy $|V_{UP} - V_{LOW}|\leq 0.2$ and $|N_{UP} - N_{LOW}|\leq 0.1$. 

In order to select the best value of the parameters,
we compare the bounds obtained by both algorithms 
after $1, 2, 5, 10, 20$, and $30$ seconds 
of computation. At each of these times, we count the number of times each 
algorithm outperforms the other by at least $1\%$ and the winner is the 
algorithm with the largest number of wins over the 6 clocked times.
It is worth noting that typically the majority of the comparisons end up as 
ties, implying that the results are not extremely sensitive to the selected
values for the parameters.

For SPARSE1, the best parameter values are $pct_{VIOL} = 0.6$ and 
$pct_{NZ} = 0.2$. For SPARSE2, they are $pct_{VIOL} = 0.6$ and 
$pct_{NZ} = 0.4$. These values are used in all experiments using
either SPARSE1 or SPARSE2 in the remainder of the paper.

\section{Computational Results}
\label{computationalresults}

In the implementation, we have used the Open Solver Interface (Osi-0.97.1) 
from COIN-OR \cite{COINOR} to create and modify the LPs 
and to interface with the LP solvers ILOG Cplex-11.1.
To compute eigenvalues and eigenvectors, we use 
the {\tt dsyevx} function provided by the LAPACK library version 3.1.1.
We also include a cut management procedure to reduce the number of 
constraints in the outer approximation LP. This procedure, applied at the 
end of each iteration, removes the cuts that are not satisfied with equality 
by the optimal solution. Note however that the constraints from the 
{\bf EXT+RLT} formulation are never removed, only constraints from
added cutting planes are possibly removed.

The machine used for the tests is a 64 bit 2.66GHz AMD processor, 
64GB of RAM memory, and Linux kernel 2.6.29.
Tolerances on the accuracy of the primal and dual solutions of the LP solver 
and LAPACK calls are set to $10^{-8}$.

The set of instances used for most experiments consists of 51 BoxQP instances 
with at most 50 variables and the 38 GLOBALLib instances as described in
Section~\ref{instances}.

For an instance $\mathcal{I}$ and a given relaxation of it, we define the 
{\it gap closed} by the relaxation as
\begin{equation}
\label{EQ:gapclosed}
100 \cdot \frac{RLT - BND}{RLT - OPT},
\end{equation}
where $BND$ and \RLT{} are the optimal value for the given relaxation and
the {\bf EXT+RLT} relaxation respectively, 
and $OPT$ is either the optimal value of $\mathcal{I}$ or the best known value
for a feasible solution. The $OPT$ values are taken from \cite{Saxena}.

\subsection{Instances}
\label{instances}
Tests are performed on a subset of instances from GL\-OB\-ALLib \cite{GLOBALLib} 
and on Box Constrained Quadratic Programs (BoxQPs) \cite{Vandenbussche}.
GLOBALLib contains 413 continuous global optimization problems of various sizes  and
types, such as BoxQPs, problems with complementarity constraints, and 
general QCQPs.
Following \cite{Saxena}, we select 160 instances from GLOBALLib having at 
most 50 variables and
that can easily be formulated as {\bf QCQP}. The conversion of a non-linear 
expression into a quadratic expression, when possible, is performed by adding 
new variables and constraints to the problem.
Additionally, bounds on the variables are derived using linear programming 
techniques and these bound are included in the formulation.
From these 160 instances in AMPL format, we substitute each bilinear term 
$x_i x_j$ by the new variable $X_{ij}$ as described for the {\bf LIFT} 
formulation. We build two collections of linearized instances in 
MPS format, one with the original precision on the coefficients and right hand 
side, and 
the second with 8-digit precision. In our experiments we used the latter.

As observed in \cite{Saxena}, using together the {\bf SDP} and {\bf RLT} 
relaxations yields stronger
bounds than those given by the {\bf RLT} relaxation only for 38 out of 
160 GLOBALLib instances. Hence, we focus on these 38 instances 
to test the effectiveness of the \SDP{} Cuts and their sparse versions.

The BoxQP collection contains 90 instances with a number of variables 
ranging from 20 to 100. Due to time limit constraints and the number of 
experiments to run, we consider only instances with a number of variables 
between 20 to 50, for a total of 51 BoxQP problems.

The converted GLOBALLib and BoxQP instances are available in MPS format from 
\cite{AQweb}.

\subsection{Effectiveness of each class of cuts}
\label{effective}

We first compare the effectiveness of the various classes of cuts when 
used in combination with the standard \SDPCUTS{}s.
For these tests, at most 1,000 cutting iterations are performed, 
at most 600 seconds are used, and operations are stopped if tailing
off is detected. More precisely, let $z_t$ be the optimal value of the linear relaxation at iteration 
$t$. The operations are halted if $t\geq 50$ and 
$z_t \ge (1-0.0001) \cdot z_{t-50}$. 
A cut purging procedure is used to remove cuts that are not tight at iteration 
$t$ if the condition $z_t \ge (1-0.0001) \cdot z_{t-1}$ is satisfied. On 
average in each iteration the algorithm generates $\frac{n^2}{2}$ cuts, 
of which only $\frac{n}{2}$ are are kept by the cut purging procedure and 
the rest are discarded.

In order to compare two different cutting plane algorithms, we compare the closed gaps 
values first after a fixed number of iterations, and second 
at several given times, for all QCQP instances at avail. Comparisons at fixed iterations indicate the quality
of the cuts, irrespective of the time used to generate them. Comparisons
at given times are useful if only limited time is available for running
the cutting plane algorithms and a good approximation of the \SDP{} cone is 
sought. The closed gaps obtained at a given point are deemed different only
if their difference is at least $g\%$ of the initial gap. We report comparisons
for $g = 1$ and $ g = 5$. Comparisons at one point is possible only if both 
algorithms reach that point. The number of problems for which this does not
happen -- because, at a given time, either result was not available 
or one of the two algorithms had already stopped, or because either algorithm 
had terminated in fewer iterations  -- is listed in the ``inc.'' (incomparable)
columns in the tables. For the remaining problems,
we report the percentage of problems for which one algorithm is better than
the other and the percentage of problems were they are tied. Finally, we also
report the average improvement in gap closed for the second algorithm 
over the first algorithm in the column labeled ``impr.''.

Tests are first performed to decide which combination 
of the SPARSE1, SPARSE2 and MINOR cuts perform best on average. Based on 
Tables \ref{bound@iter_sdpcuts_sparsecuts1_minorsdpcutsVSsdpcuts_sparsecuts2_minorsdpcuts} and
 \ref{bound@time_sdpcuts_sparsecuts1_minorsdpcutsVSsdpcuts_sparsecuts2_minorsdpcuts}
below, we conclude 
that using MINOR is useful both in terms of iteration and time,
and that the algorithm using \SDPCUTS+SPARSE2+MINOR (abbreviated {\it S2M} 
in the remainder) dominates the algorithm using 
\SDPCUTS+SPARSE1+MINOR (abbreviated {\it S1M}) both in terms of iteration
and time. Table~\ref{bound@iter_sdpcuts_sparsecuts1_minorsdpcutsVSsdpcuts_sparsecuts2_minorsdpcuts}
gives the comparison between S1M and S2M at different
iterations.
S2M dominates clearly S1M
in the very first iteration and after 200 iterations, while after the
first few iterations S1M also manages to obtain good bounds.
Table~\ref{bound@time_sdpcuts_sparsecuts1_minorsdpcutsVSsdpcuts_sparsecuts2_minorsdpcuts}
gives the comparison between these two algorithms at different
times. For comparisons with $g=1$, S1M is better than
S2M only in at most 2.25\% of the problems, while the converse
varies between roughly 50\% (at early times) and 8\% (for late
times). For $g=5$, S2M still dominates S1M in most cases.

%

%

Sparse cuts yield better bounds than using solely the standard 
\SDP{} cuts. The observed 
improvement is around 3\% and 5\% respectively for SPARSE1 and SPARSE2. 
When we are using the MINOR cuts, this value gets to 6\% and 8\% respectively 
for each type of sparsification algorithm used.
Table \ref{bound@iter_sdpcutsVSsdpcuts_sparsecuts2_minorsdpcuts}
compares \SDPCUTS{} (abbreviated by {\it S}) with S2M. The table shows 
that the sparse cuts generated by the sparsification procedures and minor 
\SDP{} cuts yield better bounds than the standard cutting plane algorithm at 
fixed iterations.
Comparisons performed at fixed times, on the other hand, show that 
considering the whole set of instances we do not get any improvement in the 
first 60 to 120 seconds of computation (see 
Table~\ref{bound@time_sdpcutsVSsdpcuts_sparsecuts2_minorsdpcuts}). Indeed S2M 
initially performs worse than the standard cutting plane algorithm, but 
after 60 to 120 seconds, it produces better bounds on average. 
In Section~\ref{appendix} detailed computational results are given in Tables 
\ref{detailedbound@iter_sdpcutsVSsdpcuts_sparsecuts2_minorsdpcuts} and
\ref{detailedbound@time_sdpcutsVSsdpcuts_sparsecuts2_minorsdpcuts} where 
for each instance we compare the duality gap closed by S and S2M at several 
iterations and times. The initial duality gap is obtained as in 
(\ref{EQ:gapclosed}) as $RLT - OPT$. We then let S2M run with no time limit
until the value $s$ obtained does not improve by at least 0.01\% over ten 
consecutive iterations. This value $s$ is
an upper bound on the value of the {\bf PSD+RLT} relaxation.
The column ``bound'' in the tables gives the
value of $RLT - s$ as a percentage of the gap $RLT - OPT$, i.e. 
an approximation of the percentage of the gap closed by the 
{\bf PSD+RLT} relaxation. The columns labeled S and S2M in the 
tables give the gap closed by the corresponding algorithms at different 
iterations.

Note that although S2M relies on numerous spectral decomposition computations,
most of its running time is spent in generating cuts and reoptimization of
the LP. For example, on the BoxQP instances with a time limit of 300 seconds, 
the average percentage of CPU time spent for obtaining spectral decompositions 
is below 21 for instances of size 30, below 15 for instances of size 40 and
below 7 for instances of size 50.

\section{Conclusions}
\label{conclusions}

\begin{figure}[h!]
\caption{Instance spar020-100-1}
\centering
\includegraphics[width=0.99\textwidth]{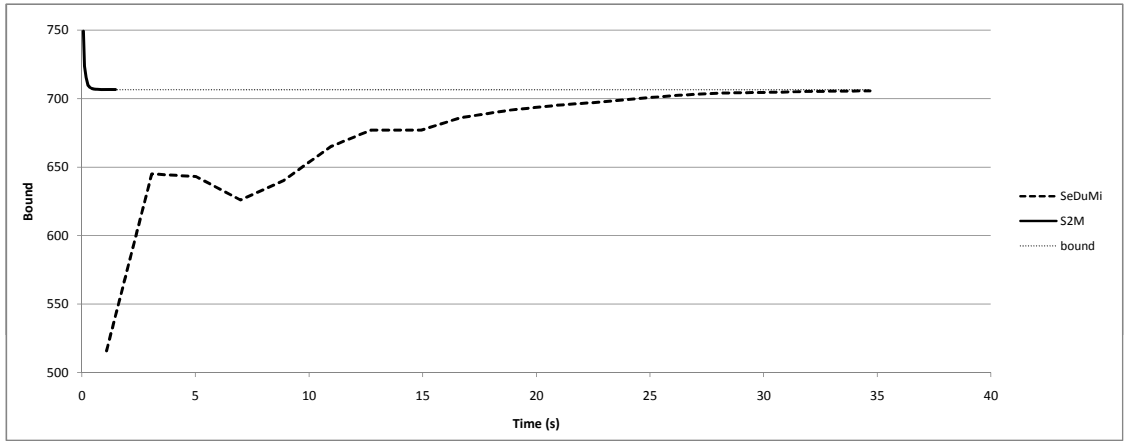}
\label{sdp0201001}
\end{figure}

\begin{figure}[h!]
\caption{Instance spar030-060-1}
\centering
\includegraphics[width=0.99\textwidth]{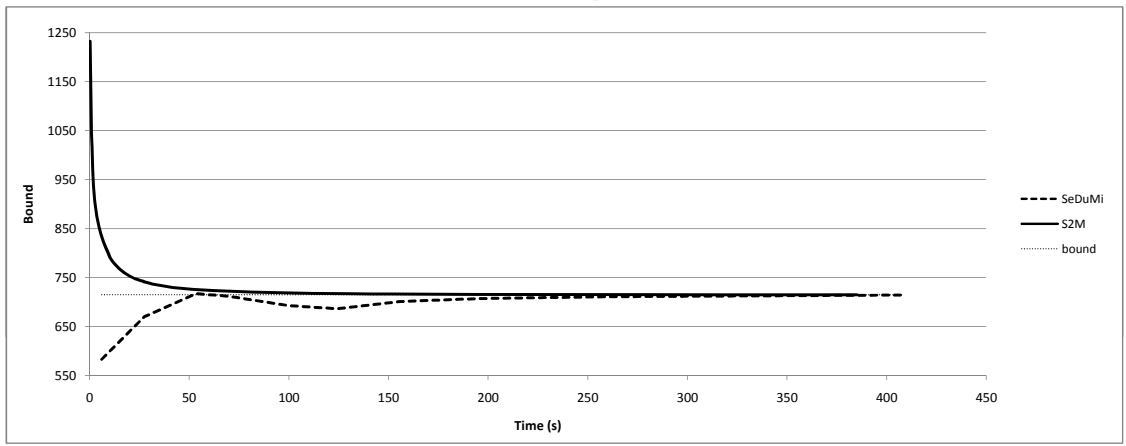}
\label{sdp0300601}
\end{figure}

This paper studies linearizations of the \SDP{} cone based on spectral
decompositions. Sparsification of eigenvectors corresponding to negative 
eigenvalues is shown to produce useful cuts in practice, in particular 
when the minor cuts are used. The goal of capturing most of the strength 
of an \SDP{} relaxation through linear inequalities is achieved, although
tailing off occurs relatively quickly. As an illustration of typical
behavior of a \SDP{} solver and our linear outer-approximation scheme,
consider the two instances, spar020-100-1 and spar030-060-1, with 
respectively 20 and 30 variables.
 We use the SDP solver {\tt SeDuMi} and S2M, keeping track at each iteration 
of the bound achieved and the time spent. 
Figure \ref{sdp0201001} and Figure \ref{sdp0300601} compare the bounds obtained 
by the two solvers at a given time.
For the small size instance spar020-100-1, we note that S2M converges to the 
bound value more than twenty times faster than {\tt SeDuMi}. In the medium 
size instance spar030-060-1 we note that S2M closes a large gap in the first 
ten to twenty iterations, and then tailing off occurs. 
To compute the exact bound, 
{\tt SeDuMi} requires 408 seconds while S2M requires 2,442 seconds to reach 
the same precision. Nevertheless, for our purposes, most of the benefits
of the \SDP{} constraints are captured in the early iterations.

Two additional improvements are possible. The first one is to use a cut
separation procedure for the RLT inequalities, avoiding their inclusion
in the initial LP and managing them as other cutting planes. This could 
potentially speed up the reoptimization of the LP. Another possibility
is to use a mix of the S and S2M algorithms, using the former in the
early iterations and then switching to the latter.

\section*{Acknowledgments}

The authors warmly thank Anureet Saxena for the useful discussions
that led to the results obtained in this work.


\newpage
\ 
\begin{table}[h!]
\caption{Comparison of S1M with S2M at several iterations.}
\label{bound@iter_sdpcuts_sparsecuts1_minorsdpcutsVSsdpcuts_sparsecuts2_minorsdpcuts}
\centering
\begin{tabular}{c|ccc|ccc|c|c}
\hline\hline
		& \multicolumn{3}{c|}{g = 1}& \multicolumn{3}{c|}{g = 5} & \ \\
\small{Iteration} 	& \small{S1M} & \small{S2M}	& \small{Tie} & \small{S1M} & \small{S2M} & \small{Tie} & \small{inc.}	& \small{impr.}\\
\hline
1	&	7.87	&	39.33	&	52.80	&	1.12	&	19.1	&	79.78	&	0.00	&	3.21\\
2	&	17.98	&	28.09	&	53.93	&	0.00	&	10.11	&	89.89	&	0.00	&	2.05\\
3	&	17.98	&	19.10	&	62.92	&	1.12	&	7.87	&	91.01	&	0.00	&	1.50\\
5	&	12.36	&	14.61	&	73.03	&	3.37	&	5.62	&	91.01	&	0.00	&	1.77\\
10	&	10.11	&	13.48	&	76.41	&	0.00	&	5.62	&	94.38	&	0.00	&	1.42\\
15	&	4.49	&	13.48	&	82.03	&	1.12	&	6.74	&	92.14	&	0.00	&	1.12\\
20	&	1.12	&	10.11	&	78.66	&	1.12	&	6.74	&	82.02	&	10.11	&	1.02\\
30	&	1.12	&	8.99	&	79.78	&	1.12	&	5.62	&	83.15	&	10.11	&	0.79\\
50	&	2.25	&	6.74	&	80.90	&	1.12	&	4.49	&	84.28	&	10.11	&	0.47\\
100	&	0.00	&	4.49	&	28.09	&	0.00	&	2.25	&	30.33	&	67.42	&	1.88\\
200	&	0.00	&	3.37	&	15.73	&	0.00	&	2.25	&	16.85	&	80.90	&	2.51\\
300	&	0.00	&	2.25	&	12.36	&	0.00	&	2.25	&	12.36	&	85.39	&	3.30\\
500	&	0.00	&	2.25	&	7.87	&	0.00	&	2.25	&	7.87	&	89.88	&	3.85\\
1000	&	0.00	&	2.25	&	3.37	&	0.00	&	2.25	&	3.37	&	94.38	&	7.43
\end{tabular}
\end{table}

\begin{table}[h!]
\caption{Comparison of S1M with S2M at several times.}
\label{bound@time_sdpcuts_sparsecuts1_minorsdpcutsVSsdpcuts_sparsecuts2_minorsdpcuts}
\centering
\begin{tabular}{c|ccc|ccc|c|c}
\hline\hline
		& \multicolumn{3}{c|}{g = 1}& \multicolumn{3}{c|}{g = 5} & \ \\
\small{Time} 	& \small{S1M} & \small{S2M}	& \small{Tie} & \small{S1M} & \small{S2M} & \small{Tie} & \small{inc.}	& \small{impr.}\\
\hline
0.5	&	3.37	&	52.81	&	12.36	&	0.00	&	43.82	&	24.72	&	31.46	&	2.77\\
1	&	0.00	&	51.68	&	14.61	&	0.00	&	40.45	&	25.84	&	33.71	&	4.35\\
2	&	0.00	&	47.19	&	15.73	&	0.00	&	39.33	&	23.59	&	37.08	&	5.89\\
3	&	1.12	&	44.94	&	14.61	&	0.00	&	34.83	&	25.84	&	39.33	&	5.11\\
5	&	1.12	&	43.82	&	15.73	&	0.00	&	38.20	&	22.47	&	39.33	&	6.07\\
10	&	1.12	&	41.58	&	16.85	&	0.00	&	24.72	&	34.83	&	40.45	&	4.97\\
15	&	2.25	&	37.08	&	16.85	&	1.12	&	21.35	&	33.71	&	43.82	&	3.64\\
20	&	1.12	&	35.96	&	16.85	&	1.12	&	17.98	&	34.83	&	46.07	&	3.49\\
30	&	1.12	&	28.09	&	22.48	&	1.12	&	16.86	&	33.71	&	48.31	&	2.99\\
60	&	1.12	&	20.23	&	28.09	&	0.00	&	12.36	&	37.08	&	50.56	&	2.62\\
120	&	0.00	&	15.73	&	32.58	&	0.00	&	10.11	&	38.20	&	51.69	&	1.73\\
180	&	0.00	&	13.49	&	32.58	&	0.00	&	5.62	&	40.45	&	53.93	&	1.19\\
300	&	0.00	&	11.24	&	31.46	&	0.00	&	3.37	&	39.33	&	57.30	&	0.92\\
600	&	0.00	&	7.86	&	24.72	&	0.00	&	0.00	&	32.58	&	67.42	&	0.72
\end{tabular}
\end{table}

\newpage
\ 
\begin{table}[h!]
\caption{Comparison of S with S2M at several iterations.}
\label{bound@iter_sdpcutsVSsdpcuts_sparsecuts2_minorsdpcuts}
\centering
\begin{tabular}{c|ccc|ccc|c|c}
\hline\hline
		& \multicolumn{3}{c|}{g = 1}& \multicolumn{3}{c|}{g = 5} & \ \\
\small{Iteration} 	& \small{S} & \small{S2M}	& \small{Tie} & \small{S} & \small{S2M} & \small{Tie} & \small{inc.}	& \small{impr.}\\
\hline
1	&	0.00	&	76.40	&	23.60	&	0.00	&	61.80	&	38.20	&	0.00	&	10.47\\
2	&	0.00	&	84.27	&	15.73	&	0.00	&	55.06	&	44.94	&	0.00	&	10.26\\
3	&	0.00	&	83.15	&	16.85	&	0.00	&	48.31	&	51.69	&	0.00	&	10.38\\
5	&	0.00	&	80.90	&	19.10	&	0.00	&	40.45	&	59.55	&	0.00	&	10.09\\
10	&	1.12	&	71.91	&	26.97	&	0.00	&	41.57	&	58.43	&	0.00	&	8.87\\
15	&	1.12	&	60.67	&	38.21	&	1.12	&	35.96	&	62.92	&	0.00	&	7.49\\
20	&	1.12	&	53.93	&	40.45	&	1.12	&	29.21	&	65.17	&	4.50	&	6.22\\
30	&	1.12	&	34.83	&	53.93	&	0.00	&	16.85	&	73.03	&	10.12	&	5.04\\
50	&	1.12	&	25.84	&	62.92	&	0.00	&	13.48	&	76.40	&	10.12	&	3.75\\
100	&	1.12	&	8.99	&	21.35	&	0.00	&	5.62	&	25.84	&	68.54	&	5.57\\
200	&	0.00	&	5.62	&	8.99	&	0.00	&	3.37	&	11.24	&	85.39	&	7.66\\
300	&	0.00	&	3.37	&	7.87	&	0.00	&	3.37	&	7.87	&	88.76	&	8.86\\
500	&	0.00	&	3.37	&	5.62	&	0.00	&	3.37	&	5.62	&	91.01	&	8.72\\
1000	&	0.00	&	2.25	&	0.00	&	0.00	&	2.25	&	0.00	&	97.75	&	26.00
\end{tabular}
\end{table}

\begin{table}[h!]
\caption{Comparison of S with S2M at several times.}
\label{bound@time_sdpcutsVSsdpcuts_sparsecuts2_minorsdpcuts}
\centering
\begin{tabular}{c|ccc|ccc|c|c}
\hline\hline
		& \multicolumn{3}{c|}{g = 1}& \multicolumn{3}{c|}{g = 5} & \ \\
\small{Time} 	& \small{S} & \small{S2M}	& \small{Tie} & \small{S} & \small{S2M} & \small{Tie} & \small{inc.}	& \small{impr.}\\
\hline
0.5	&	41.57	&	17.98	&	5.62	&	41.57	&	17.98	&	5.62	&	34.83	&	-9.42\\
1	&	41.57	&	14.61	&	5.62	&	39.33	&	13.48	&	8.99	&	38.20	&	-8.66\\
2	&	42.70	&	10.11	&	6.74	&	29.21	&	8.99	&	21.35	&	40.45	&	-8.73\\
3	&	41.57	&	8.99	&	8.99	&	31.46	&	6.74	&	21.35	&	40.45	&	-8.78\\
5	&	35.96	&	7.87	&	15.72	&	33.71	&	5.62	&	20.22	&	40.45	&	-7.87\\
10	&	34.84	&	7.87	&	13.48	&	30.34	&	4.50	&	21.35	&	43.81	&	-5.95\\
15	&	37.07	&	5.62	&	11.24	&	22.47	&	2.25	&	29.21	&	46.07	&	-5.48\\
20	&	37.07	&	5.62	&	8.99	&	17.98	&	1.12	&	32.58	&	48.32	&	-4.99\\
30	&	30.34	&	5.62	&	15.72	&	11.24	&	1.12	&	39.32	&	48.32	&	-3.9\\
60	&	11.24	&	12.36	&	25.84	&	11.24	&	2.25	&	35.95	&	50.56	&	-1.15\\
120	&	8.99	&	12.36	&	24.72	&	2.25	&	2.25	&	41.57	&	53.93	&	 0.48\\
180	&	2.25	&	14.61	&	29.21	&	0.00	&	4.50	&	41.57	&	53.93	&	 1.09\\
300	&	0.00	&	15.73	&	26.97	&	0.00	&	6.74	&	35.96	&	57.30	&	 1.60\\
600	&	0.00	&	14.61	&	13.48	&	0.00	&	5.62	&	22.47	&	71.91	&	 2.73
\end{tabular}
\end{table}

\begin{landscape}
 \section{Appendix}
 \label{appendix}

\scriptsize
\begin{center}
\begin{table}[h!] 
\caption{Duality gap closed at several iterations for each instance.}
\label{detailedbound@iter_sdpcutsVSsdpcuts_sparsecuts2_minorsdpcuts}
\end{table}
\begin{longtable}{r|rr|r|rr|rr|rr}

  \hline\hline\\[-2ex]
	&	&	&	&	\multicolumn{2}{|c|}{iter. 2}			&	\multicolumn{2}{|c|}{iter. 10}			&	\multicolumn{2}{|c}{iter. 50} \\
Instance	&	$|x|$	&	$|y|$	&	bound	&	S	&	S2M	&	S	&	S2M	&	S	&	S2M	\\ 
\hline\\[-1.8ex]
\endfirsthead 

\multicolumn{10}{c}{{\tablename} \ref{detailedbound@iter_sdpcutsVSsdpcuts_sparsecuts2_minorsdpcuts} -- Continued} \\[0.5ex]
  \hline\hline\\[-2ex]
	&	&	&	&	\multicolumn{2}{|c|}{iter. 2}			&	\multicolumn{2}{|c|}{iter. 10}			&	\multicolumn{2}{|c}{iter. 50} \\
Instance	&	$|x|$	&	$|y|$	&	bound	&	S	&	S2M	&	S	&	S2M	&	S	&	S2M	\\
\hline\\[-1.8ex]
\endhead

\multicolumn{10}{l}{{Continued on Next Page\ldots}} \\
\endfoot

\multicolumn{10}{c}{}\\[-1.8ex]
\endlastfoot

circle	&	3	&	0	&	45.79	&	0.00	&	0.00	&	10.97	&	41.31	&	45.77	&	45.79	\\ 
dispatch	&	3	&	1	&	100.00	&	25.59	&	27.92	&	37.25	&	35.76	&	95.90	&	92.17	\\ 
ex2\_1\_10	&	20	&	0	&	22.05	&	3.93	&	8.65	&	15.93	&	21.05	&	22.05	&	22.05	\\ 
ex3\_1\_2	&	5	&	0	&	49.75	&	49.75	&	49.75	&	49.75	&	49.75	&	49.75	&	49.75	\\ 
ex4\_1\_1	&	3	&	0	&	100.00	&	99.81	&	99.84	&	100.00	&	100.00	&	100.00	&	100.00	\\ 
ex4\_1\_3	&	3	&	0	&	56.40	&	0.00	&	0.00	&	51.19	&	51.19	&	56.40	&	56.40	\\ 
ex4\_1\_4	&	3	&	0	&	100.00	&	22.33	&	42.78	&	98.98	&	99.98	&	100.00	&	100.00	\\ 
ex4\_1\_6	&	3	&	0	&	100.00	&	69.44	&	69.87	&	92.62	&	99.94	&	100.00	&	100.00	\\ 
ex4\_1\_7	&	3	&	0	&	100.00	&	18.00	&	48.17	&	96.86	&	99.90	&	100.00	&	100.00	\\ 
ex4\_1\_8	&	3	&	0	&	100.00	&	56.90	&	81.93	&	99.76	&	99.93	&	100.00	&	100.00	\\ 
ex8\_1\_4	&	4	&	0	&	100.00	&	94.91	&	95.19	&	99.98	&	100.00	&	100.00	&	100.00	\\ 
ex8\_1\_5	&	5	&	0	&	68.26	&	32.32	&	39.17	&	59.01	&	66.76	&	68.00	&	68.25	\\ 
ex8\_1\_7	&	9	&	0	&	77.43	&	3.04	&	33.75	&	33.13	&	53.44	&	64.03	&	75.38	\\ 
ex8\_4\_1	&	21	&	1	&	91.81	&	4.45	&	21.80	&	18.60	&	45.08	&	38.07	&	69.83	\\ 
ex9\_2\_1	&	10	&	0	&	54.52	&	0.00	&	42.55	&	0.01	&	50.13	&	0.01	&	51.90	\\ 
ex9\_2\_2	&	10	&	0	&	70.37	&	0.00	&	14.08	&	2.34	&	51.97	&	7.12	&	69.41	\\ 
ex9\_2\_4	&	6	&	2	&	99.87	&	0.00	&	24.84	&	25.24	&	99.85	&	86.37	&	99.87	\\ 
ex9\_2\_6	&	16	&	0	&	99.88	&	3.50	&	99.42	&	23.09	&	99.86	&	62.32	&	99.88	\\ 
ex9\_2\_7	&	10	&	0	&	42.30	&	0.00	&	4.59	&	0.00	&	27.34	&	3.14	&	34.91	\\ 
himmel11	&	5	&	4	&	49.75	&	49.75	&	49.75	&	49.75	&	49.75	&	49.75	&	49.75	\\ 
hydro	&	12	&	19	&	52.06	&	0.00	&	20.87	&	21.95	&	29.03	&	26.04	&	31.39	\\ 
mathopt1	&	4	&	0	&	100.00	&	95.76	&	100.00	&	99.96	&	100.00	&	100.00	&	100.00	\\ 
mathopt2	&	3	&	0	&	100.00	&	99.84	&	99.93	&	100.00	&	100.00	&	100.00	&	100.00	\\ 
meanvar	&	7	&	1	&	100.00	&	0.00	&	0.00	&	78.35	&	95.84	&	100.00	&	100.00	\\ 
nemhaus	&	5	&	0	&	53.97	&	26.00	&	26.41	&	48.49	&	50.16	&	53.87	&	53.96	\\ 
prob06	&	2	&	0	&	100.00	&	90.61	&	92.39	&	98.39	&	98.39	&	98.39	&	98.39	\\ 
prob09	&	3	&	1	&	100.00	&	0.00	&	99.00	&	61.14	&	99.96	&	99.64	&	100.00	\\ 
process	&	9	&	3	&	8.00	&	0.00	&	4.25	&	0.00	&	4.98	&	0.00	&	5.73	\\ 
qp1	&	50	&	0	&	100.00	&	79.59	&	89.09	&	93.89	&	99.77	&	98.93	&	100.00	\\ 
qp2	&	50	&	0	&	100.00	&	55.94	&	70.99	&	82.42	&	93.92	&	93.04	&	99.35	\\ 
rbrock	&	3	&	0	&	100.00	&	97.48	&	100.00	&	99.96	&	100.00	&	100.00	&	100.00	\\ 
st\_e10	&	3	&	1	&	100.00	&	56.90	&	81.93	&	99.76	&	99.93	&	100.00	&	100.00	\\ 
st\_e18	&	2	&	0	&	100.00	&	0.00	&	0.00	&	98.72	&	98.72	&	100.00	&	100.00	\\ 
st\_e19	&	3	&	1	&	93.51	&	5.14	&	15.93	&	29.97	&	60.10	&	93.40	&	93.50	\\ 
st\_e25	&	4	&	0	&	87.55	&	55.80	&	55.80	&	87.02	&	87.01	&	87.23	&	87.23	\\ 
st\_e28	&	5	&	4	&	49.75	&	49.75	&	49.75	&	49.75	&	49.75	&	49.75	&	49.75	\\ 
st\_iqpbk1	&	8	&	0	&	97.99	&	71.99	&	76.69	&	97.20	&	97.95	&	97.99	&	97.99	\\ 
st\_iqpbk2	&	8	&	0	&	97.93	&	70.55	&	75.16	&	94.93	&	97.52	&	97.93	&	97.93	\\ 
spar020-100-1	&	20	&	0	&	100.00	&	91.15	&	94.64	&	99.77	&	99.99	&	100.00	&	100.00	\\ 
spar020-100-2	&	20	&	0	&	99.70	&	90.12	&	92.64	&	98.17	&	99.32	&	99.66	&	99.69	\\ 
spar020-100-3	&	20	&	0	&	100.00	&	96.96	&	98.51	&	100.00	&	100.00	&	100.00	&	100.00	\\ 
spar030-060-1	&	30	&	0	&	98.87	&	43.53	&	53.64	&	79.61	&	87.39	&	93.90	&	97.14	\\ 
spar030-060-2	&	30	&	0	&	100.00	&	80.74	&	89.73	&	99.89	&	100.00	&	100.00	&	100.00	\\ 
spar030-060-3	&	30	&	0	&	99.40	&	67.43	&	71.94	&	91.48	&	95.68	&	98.75	&	99.26	\\ 
spar030-070-1	&	30	&	0	&	97.99	&	49.05	&	54.94	&	76.54	&	86.51	&	91.15	&	95.68	\\ 
spar030-070-2	&	30	&	0	&	100.00	&	81.19	&	85.82	&	99.26	&	99.99	&	100.00	&	100.00	\\ 
spar030-070-3	&	30	&	0	&	99.98	&	85.97	&	87.43	&	98.44	&	99.52	&	99.92	&	99.97	\\ 
spar030-080-1	&	30	&	0	&	98.99	&	64.44	&	70.99	&	87.32	&	92.11	&	96.23	&	98.01	\\ 
spar030-080-2	&	30	&	0	&	100.00	&	92.78	&	95.45	&	100.00	&	100.00	&	100.00	&	100.00	\\ 
spar030-080-3	&	30	&	0	&	100.00	&	92.71	&	94.18	&	99.99	&	100.00	&	100.00	&	100.00	\\ 
spar030-090-1	&	30	&	0	&	100.00	&	80.37	&	86.35	&	97.27	&	99.30	&	100.00	&	100.00	\\ 
spar030-090-2	&	30	&	0	&	100.00	&	86.09	&	89.26	&	98.13	&	99.65	&	100.00	&	100.00	\\ 
spar030-090-3	&	30	&	0	&	100.00	&	90.65	&	91.56	&	99.97	&	100.00	&	100.00	&	100.00	\\ 
spar030-100-1	&	30	&	0	&	100.00	&	77.28	&	83.25	&	95.20	&	98.30	&	99.85	&	100.00	\\ 
spar030-100-2	&	30	&	0	&	99.96	&	76.78	&	81.65	&	93.44	&	96.84	&	98.70	&	99.72	\\ 
spar030-100-3	&	30	&	0	&	99.85	&	86.82	&	88.74	&	97.45	&	98.75	&	99.75	&	99.83	\\ 
spar040-030-1	&	40	&	0	&	100.00	&	25.60	&	41.96	&	73.59	&	84.72	&	99.13	&	100.00	\\ 
spar040-030-2	&	40	&	0	&	100.00	&	30.93	&	53.39	&	79.34	&	95.62	&	99.46	&	100.00	\\ 
spar040-030-3	&	40	&	0	&	100.00	&	9.21	&	31.38	&	66.46	&	86.62	&	98.53	&	100.00	\\ 
spar040-040-1	&	40	&	0	&	96.74	&	23.62	&	29.03	&	63.04	&	75.93	&	85.93	&	93.29	\\ 
spar040-040-2	&	40	&	0	&	100.00	&	33.17	&	48.87	&	89.08	&	97.94	&	100.00	&	100.00	\\ 
spar040-040-3	&	40	&	0	&	99.18	&	21.77	&	30.31	&	70.44	&	80.96	&	91.37	&	96.69	\\ 
spar040-050-1	&	40	&	0	&	99.42	&	35.62	&	44.87	&	73.11	&	84.05	&	92.81	&	97.21	\\ 
spar040-050-2	&	40	&	0	&	99.48	&	36.79	&	47.68	&	82.38	&	91.27	&	97.26	&	98.93	\\ 
spar040-050-3	&	40	&	0	&	100.00	&	41.91	&	51.72	&	84.04	&	90.70	&	96.88	&	99.34	\\ 
spar040-060-1	&	40	&	0	&	98.09	&	46.22	&	52.89	&	81.65	&	87.28	&	92.39	&	95.97	\\ 
spar040-060-2	&	40	&	0	&	100.00	&	63.02	&	72.87	&	94.09	&	97.66	&	99.78	&	100.00	\\ 
spar040-060-3	&	40	&	0	&	100.00	&	78.09	&	87.91	&	99.30	&	99.99	&	100.00	&	100.00	\\ 
spar040-070-1	&	40	&	0	&	100.00	&	64.02	&	71.33	&	93.92	&	97.35	&	99.77	&	100.00	\\ 
spar040-070-2	&	40	&	0	&	100.00	&	67.49	&	76.78	&	95.12	&	97.97	&	99.97	&	100.00	\\ 
spar040-070-3	&	40	&	0	&	100.00	&	70.13	&	79.43	&	95.65	&	97.99	&	99.75	&	100.00	\\ 
spar040-080-1	&	40	&	0	&	100.00	&	63.06	&	69.40	&	91.09	&	95.44	&	99.00	&	99.97	\\ 
spar040-080-2	&	40	&	0	&	100.00	&	71.42	&	79.77	&	94.98	&	97.62	&	99.92	&	100.00	\\ 
spar040-080-3	&	40	&	0	&	99.99	&	83.93	&	88.65	&	97.76	&	98.86	&	99.81	&	99.95	\\ 
spar040-090-1	&	40	&	0	&	100.00	&	75.73	&	79.96	&	95.34	&	97.43	&	99.46	&	99.91	\\ 
spar040-090-2	&	40	&	0	&	99.97	&	76.39	&	80.97	&	95.16	&	96.72	&	99.20	&	99.81	\\ 
spar040-090-3	&	40	&	0	&	100.00	&	84.90	&	87.04	&	98.33	&	99.52	&	100.00	&	100.00	\\ 
spar040-100-1	&	40	&	0	&	100.00	&	87.64	&	90.43	&	98.27	&	99.35	&	99.98	&	100.00	\\ 
spar040-100-2	&	40	&	0	&	99.87	&	79.78	&	83.02	&	94.58	&	96.76	&	98.74	&	99.50	\\ 
spar040-100-3	&	40	&	0	&	98.70	&	72.69	&	78.31	&	90.83	&	93.03	&	95.84	&	97.36	\\ 
spar050-030-1	&	50	&	0	&	100.00	&	3.11	&	17.60	&	58.23	&	79.98	&	-	&	-	\\ 
spar050-030-2	&	50	&	0	&	99.27	&	1.35	&	16.67	&	51.11	&	70.58	&	-	&	-	\\ 
spar050-030-3	&	50	&	0	&	99.29	&	0.08	&	13.63	&	50.19	&	67.46	&	-	&	-	\\ 
spar050-040-1	&	50	&	0	&	100.00	&	23.13	&	30.86	&	72.10	&	81.73	&	-	&	-	\\ 
spar050-040-2	&	50	&	0	&	99.39	&	21.89	&	34.45	&	71.24	&	81.63	&	-	&	-	\\ 
spar050-040-3	&	50	&	0	&	100.00	&	27.18	&	37.42	&	83.96	&	91.70	&	-	&	-	\\ 
spar050-050-1	&	50	&	0	&	93.02	&	25.24	&	33.77	&	61.42	&	68.75	&	-	&	-	\\ 
spar050-050-2	&	50	&	0	&	98.74	&	32.10	&	41.26	&	77.48	&	83.48	&	-	&	-	\\ 
spar050-050-3	&	50	&	0	&	98.84	&	38.57	&	44.67	&	80.97	&	85.36	&	-	&	-	\\ 
\hline
Average	&	-	&	-	&	-	&	48.75	&	59.00	&	75.53	&	84.39	&	85.85	&	89.60
\end{longtable}
\addtocounter{table}{-1} 
\end{center}
\normalsize


\scriptsize
\begin{center}
\begin{table} 
\caption{Duality gap closed at several times for each instance. (Instances solved in less than 1 second are not shown)}
\label{detailedbound@time_sdpcutsVSsdpcuts_sparsecuts2_minorsdpcuts}
\end{table}
\begin{longtable}{r|r|rr|rr|rr|rr|rr}

  \hline\hline\\[-2ex]
   		&	& \multicolumn{2}{|c|}{1 s}	& \multicolumn{2}{|c|}{60 s}	& \multicolumn{2}{|c|}{180 s}	& \multicolumn{2}{|c}{300 s} & \multicolumn{2}{|c}{600 s}  \\
  Instance	&	bound	&	S	&	S2M	&	S	&	S2M	&	S	&	S2M	&	S	&	S2M	&	S	&	S2M\\ \hline\\[-1.8ex]
\endfirsthead 

\multicolumn{10}{c}{{\tablename} \ref{detailedbound@time_sdpcutsVSsdpcuts_sparsecuts2_minorsdpcuts} -- Continued} \\[0.5ex]
  \hline\hline\\[-2ex]
   		&	& \multicolumn{2}{|c|}{1 s}	& \multicolumn{2}{|c|}{60 s}	& \multicolumn{2}{|c|}{180 s}	& \multicolumn{2}{|c}{300 s} & \multicolumn{2}{|c}{600 s}  \\
  Instance	&	bound	&	S	&	S2M	&	S	&	S2M	&	S	&	S2M	&	S	&	S2M	&	S	&	S2M\\ \hline\\[-1.8ex]
\endhead

\multicolumn{10}{l}{{Continued on Next Page\ldots}} \\
\endfoot

\multicolumn{12}{c}{}\\[-1.8ex] 
\endlastfoot

ex4\_1\_4	&	100.00	&	-	&	100.00	&	-	&	-	&	-	&	-	&	-	&	-	&	-	&	-	 \\ 
ex8\_1\_4	&	100.00	&	-	&	100.00	&	-	&	-	&	-	&	-	&	-	&	-	&	-	&	-	 \\ 
ex8\_1\_7	&	77.43	&	77.43	&	77.37	&	-	&	-	&	-	&	-	&	-	&	-	&	-	&	-	 \\ 
ex8\_4\_1	&	91.81	&	28.14	&	36.24	&	61.60	&	90.43	&	-	&	-	&	-	&	-	&	-	&	-	 \\ 
ex9\_2\_2	&	70.37	&	-	&	70.35	&	-	&	-	&	-	&	-	&	-	&	-	&	-	&	-	 \\ 
ex9\_2\_6	&	99.88	&	96.28	&	-	&	-	&	-	&	-	&	-	&	-	&	-	&	-	&	-	 \\ 
hydro	&	52.06	&	26.43	&	31.46	&	-	&	-	&	-	&	-	&	-	&	-	&	-	&	-	 \\ 
mathopt2	&	100.00	&	-	&	100.00	&	-	&	-	&	-	&	-	&	-	&	-	&	-	&	-	 \\ 
process	&	8.00	&	-	&	7.66	&	-	&	-	&	-	&	-	&	-	&	-	&	-	&	-	 \\ 
qp1	&	100.00	&	79.99	&	80.28	&	98.22	&	99.52	&	99.73	&	99.96	&	99.92	&	99.98	&	99.99	&	100.00	 \\ 
qp2	&	100.00	&	55.82	&	55.27	&	91.74	&	95.56	&	95.86	&	98.69	&	97.41	&	99.66	&	98.80	&	100.00	 \\ 
spar020-100-1	&	100.00	&	100.00	&	100.00	&	-	&	-	&	-	&	-	&	-	&	-	&	-	&	-	 \\ 
spar020-100-2	&	99.70	&	99.67	&	99.61	&	-	&	-	&	-	&	-	&	-	&	-	&	-	&	-	 \\ 
spar020-100-3	&	100.00	&	-	&	100.00	&	-	&	-	&	-	&	-	&	-	&	-	&	-	&	-	 \\ 
spar030-060-1	&	98.87	&	69.98	&	58.72	&	96.53	&	97.61	&	98.45	&	98.70	&	98.68	&	98.82	&	-	&	-	 \\ 
spar030-060-2	&	100.00	&	96.52	&	91.05	&	-	&	-	&	-	&	-	&	-	&	-	&	-	&	-	 \\ 
spar030-060-3	&	99.40	&	82.99	&	76.15	&	99.27	&	99.32	&	99.38	&	99.39	&	99.39	&	99.40	&	99.40	&	99.40	 \\ 
spar030-070-1	&	97.99	&	69.81	&	60.36	&	94.50	&	96.38	&	97.29	&	97.73	&	97.70	&	97.91	&	-	&	97.98	 \\ 
spar030-070-2	&	100.00	&	96.05	&	87.93	&	-	&	-	&	-	&	-	&	-	&	-	&	-	&	-	 \\ 
spar030-070-3	&	99.98	&	96.26	&	90.42	&	99.98	&	99.98	&	99.98	&	99.98	&	-	&	99.98	&	-	&	-	 \\ 
spar030-080-1	&	98.99	&	83.36	&	74.42	&	97.80	&	98.11	&	98.74	&	98.88	&	98.89	&	98.96	&	-	&	98.99	 \\ 
spar030-080-2	&	100.00	&	99.83	&	96.70	&	-	&	-	&	-	&	-	&	-	&	-	&	-	&	-	 \\ 
spar030-080-3	&	100.00	&	99.88	&	95.87	&	-	&	-	&	-	&	-	&	-	&	-	&	-	&	-	 \\ 
spar030-090-1	&	100.00	&	92.86	&	87.69	&	-	&	-	&	-	&	-	&	-	&	-	&	-	&	-	 \\ 
spar030-090-2	&	100.00	&	93.80	&	88.46	&	-	&	100.00	&	-	&	-	&	-	&	-	&	-	&	-	 \\ 
spar030-090-3	&	100.00	&	97.78	&	91.35	&	-	&	-	&	-	&	-	&	-	&	-	&	-	&	-	 \\ 
spar030-100-1	&	100.00	&	91.04	&	84.34	&	100.00	&	100.00	&	-	&	-	&	-	&	-	&	-	&	-	 \\ 
spar030-100-2	&	99.96	&	90.21	&	83.14	&	99.56	&	99.75	&	99.91	&	99.95	&	99.95	&	99.96	&	-	&	99.96	 \\ 
spar030-100-3	&	99.85	&	94.26	&	89.55	&	99.84	&	99.84	&	99.85	&	99.85	&	99.85	&	99.85	&	99.85	&	99.85	 \\ 
spar040-030-1	&	100.00	&	28.97	&	40.51	&	89.30	&	84.19	&	99.06	&	99.98	&	99.98	&	100.00	&	-	&	100.00	 \\ 
spar040-030-2	&	100.00	&	31.97	&	48.01	&	94.01	&	96.39	&	99.58	&	99.98	&	99.99	&	100.00	&	-	&	-	 \\ 
spar040-030-3	&	100.00	&	9.20	&	27.59	&	81.66	&	85.43	&	97.25	&	99.86	&	99.81	&	100.00	&	100.00	&	-	 \\ 
spar040-040-1	&	96.74	&	19.38	&	22.90	&	70.35	&	75.45	&	80.73	&	88.63	&	85.34	&	92.29	&	90.79	&	94.74	 \\ 
spar040-040-2	&	100.00	&	24.51	&	29.87	&	98.63	&	98.60	&	100.00	&	100.00	&	-	&	-	&	-	&	-	 \\ 
spar040-040-3	&	99.18	&	20.88	&	21.31	&	78.28	&	79.31	&	86.02	&	91.22	&	89.52	&	95.04	&	94.07	&	97.71	 \\ 
spar040-050-1	&	99.42	&	28.96	&	21.27	&	80.18	&	84.01	&	88.70	&	94.62	&	92.75	&	96.71	&	96.53	&	98.32	 \\ 
spar040-050-2	&	99.48	&	29.52	&	16.91	&	91.33	&	91.42	&	97.01	&	97.97	&	98.26	&	98.87	&	-	&	99.31	 \\ 
spar040-050-3	&	100.00	&	28.67	&	19.81	&	90.03	&	90.72	&	95.68	&	97.51	&	97.49	&	99.08	&	98.92	&	99.89	 \\ 
spar040-060-1	&	98.09	&	37.16	&	17.10	&	86.26	&	87.13	&	90.18	&	93.50	&	92.25	&	95.32	&	95.05	&	96.84	 \\ 
spar040-060-2	&	100.00	&	39.57	&	22.83	&	98.09	&	98.22	&	99.90	&	99.96	&	100.00	&	100.00	&	100.00	&	-	 \\ 
spar040-060-3	&	100.00	&	52.41	&	30.57	&	100.00	&	99.99	&	-	&	-	&	-	&	-	&	-	&	-	 \\ 
spar040-070-1	&	100.00	&	50.01	&	21.79	&	97.74	&	97.78	&	99.80	&	99.87	&	99.97	&	99.99	&	100.00	&	100.00	 \\ 
spar040-070-2	&	100.00	&	47.57	&	25.19	&	98.81	&	98.46	&	99.99	&	99.99	&	100.00	&	100.00	&	-	&	-	 \\ 
spar040-070-3	&	100.00	&	47.22	&	21.95	&	98.96	&	98.70	&	99.88	&	99.92	&	99.98	&	100.00	&	100.00	&	100.00	 \\ 
spar040-080-1	&	100.00	&	51.66	&	28.00	&	95.13	&	95.38	&	98.29	&	99.05	&	99.09	&	99.74	&	99.77	&	99.99	 \\ 
spar040-080-2	&	100.00	&	52.24	&	25.94	&	98.71	&	98.31	&	99.95	&	99.97	&	100.00	&	100.00	&	-	&	-	 \\ 
spar040-080-3	&	99.99	&	56.05	&	26.98	&	99.54	&	99.25	&	99.89	&	99.88	&	99.94	&	99.95	&	99.97	&	99.98	 \\ 
spar040-090-1	&	100.00	&	59.71	&	28.17	&	98.10	&	97.86	&	99.43	&	99.61	&	99.70	&	99.86	&	99.90	&	99.99	 \\ 
spar040-090-2	&	99.97	&	59.14	&	29.82	&	97.83	&	97.70	&	99.34	&	99.58	&	99.68	&	99.81	&	99.86	&	99.93	 \\ 
spar040-090-3	&	100.00	&	63.07	&	34.62	&	99.94	&	99.85	&	100.00	&	100.00	&	-	&	-	&	-	&	-	 \\ 
spar040-100-1	&	100.00	&	69.47	&	28.24	&	99.66	&	99.47	&	99.99	&	99.99	&	100.00	&	100.00	&	-	&	-	 \\ 
spar040-100-2	&	99.87	&	65.27	&	26.07	&	97.34	&	96.87	&	98.60	&	98.98	&	99.02	&	99.39	&	99.44	&	99.69	 \\ 
spar040-100-3	&	98.70	&	61.40	&	29.61	&	93.01	&	93.17	&	94.91	&	96.02	&	95.81	&	97.00	&	96.84	&	97.77	 \\ 
spar050-030-1	&	100.00	&	0.37	&	3.63	&	54.46	&	37.52	&	70.10	&	73.34	&	76.87	&	84.75	&	86.23	&	96.33	 \\ 
spar050-030-2	&	99.27	&	0.08	&	2.79	&	44.68	&	38.62	&	59.58	&	64.94	&	67.79	&	74.98	&	77.02	&	86.58	 \\ 
spar050-030-3	&	99.29	&	0.00	&	2.75	&	44.32	&	32.31	&	57.13	&	59.07	&	62.54	&	68.99	&	71.18	&	82.86	 \\ 
spar050-040-1	&	100.00	&	3.76	&	1.77	&	69.97	&	56.87	&	77.15	&	78.30	&	80.31	&	84.30	&	84.90	&	91.79	 \\ 
spar050-040-2	&	99.39	&	2.08	&	2.84	&	68.64	&	58.47	&	77.72	&	77.61	&	81.54	&	83.63	&	86.40	&	90.94	 \\ 
spar050-040-3	&	100.00	&	1.76	&	2.31	&	79.44	&	65.71	&	89.73	&	87.74	&	92.67	&	93.00	&	95.99	&	97.69	 \\ 
spar050-050-1	&	93.02	&	4.91	&	1.84	&	60.64	&	53.28	&	65.52	&	66.42	&	66.81	&	70.38	&	68.45	&	74.76	 \\ 
spar050-050-2	&	98.74	&	6.18	&	3.39	&	76.56	&	68.33	&	82.34	&	82.21	&	84.94	&	86.52	&	-	&	91.34	 \\ 
spar050-050-3	&	98.84	&	6.12	&	2.82	&	79.38	&	69.23	&	84.95	&	83.23	&	86.99	&	86.98	&	89.77	&	91.57	 \\  \hline
Average	&	-	&	51.45	&	42.96	&	87.50	&	86.38	&	92.14	&	93.22	&	93.18	&	94.77	&	93.16	&	95.86\\
\end{longtable}
\addtocounter{table}{-1} 
\end{center}
\normalsize
\end{landscape}


\begin{thebibliography}{99}


\bibitem{al-khayyal}
F. A. Al-Khayyal and J. E. Falk, \newblock Jointly constrained biconvex programming. {\em Math. Oper. Res.} 8, pp.~273-286, 1983.

\bibitem{Anstreicher}
 K. M. Anstreicher, \newblock Semidefinite Programming versus the Reformulation Linearization Technique for Nonconvex Quadratically Constrained Quadratic Programming. Pre-print. {\em Optimization Online}, May 2007.
Available at 
\href{http://www.optimization-online.org/DB_HTML/2007/05/1655..html}
     {\tt\small{http://www.optimization-online.org/DB\_HTML/2007/05/1655.html}}


\bibitem{Balas}
E. Balas, \newblock Disjunctive programming: properties of the convex hull of feasible points. {\em Disc. Appl. Math.} 89, 1998.

\bibitem{nlp}
M. S. Bazaraa, H. D. Sherali and C. M. Shetty, \newblock Nonlinear Programming: Theory and Algorithms. Wiley, 2006.

\bibitem{BoydVandenberghe}
S. Boyd, L. Vandenberghe, \newblock Convex Optimization. Cambridge University Press, 2004.

\bibitem{BurerLetchford}
S. Burer and A. Letchford, \newblock On Non-Convex Quadratic Programming with Box Constraints. {\em Optimization Online}, July 2008.
Available at\\
\href{http://www.optimization-online.org/DB_HTML/2008/07/2030.html}
     {\tt\small{http://www.optimization-online.org/DB\_HTML/2008/07/2030.html}}

\bibitem{CSDP} 
B. Borchers, \newblock CSDP, A C Library for Semidefinite Programming, {\em Optimization Methods and
Software} 11(1), pp.~613-623, 1999.

\bibitem{COINOR}
COmputational INfrastructure for Operations Research (COIN-OR).\\ \newblock \href{http://www.coin-or.org}{\tt\small{http://www.coin-or.org}}

\bibitem{DSDP}
S. J. Benson and Y. Ye, \newblock {DSDP5}: Software For Semidefinite Programming. Available at \href{http://www-unix.mcs.anl.gov/DSDP}{\tt{\small{http://www-unix.mcs.anl.gov/DSDP}}}

\bibitem{GLOBALLib} 
GamsWorld Global Optimization library.\\\href{http://www.gamsworld.org/global/globallib/globalstat.htm}{\tt\small{http://www.gamsworld.org/global/globallib/globalstat.htm}}

\bibitem{Lovasz}
L. Lov\'{a}sz and A. Schrijver, \newblock Cones of matrices and set-functions and 0-1 optimization. {\em SIAM Journal on Optimization}, May 1991

\bibitem{mccormick}
G.P. McCormick, \newblock Nonlinear programming: theory, algorithms and applications. John Wiley \& sons, 1983.

\bibitem{AQweb}
  \href{http://www.andrew.cmu.edu/user/aqualizz/research/MIQCP}
       {\tt\small{http://www.andrew.cmu.edu/user/aqualizz/research/MIQCP}}



\bibitem{Saxena}
A. Saxena, P. Bonami and J. Lee, Convex Relaxations of Non-Convex Mixed Integer Quadratically Constrained Programs: Extended Formulations, 
2009. To appear in {\em Mathematical Programming}.


\bibitem{anureet}
\bysame,
 Convex Relaxations
  of Non-Convex Mixed Integer Quadratically Constrained Programs:
  Projected Formulations, Optimization Online, November
  2008. Available at\\
  \href{http://www.optimization-online.org/DB_HTML/2008/11/2145.html}
       {\tt{\small{http://www.optimization-online.org/DB\_HTML/2008/11/2145.html}}}

\bibitem{SeDuMi}
J. F. Sturm, \newblock SeDuMi: An Optimization Package over Symmetric
Cones. Available at 
\href{http://sedumi.mcmaster.ca}
     {\tt\small{http://sedumi.mcmaster.ca}}

\bibitem{SheraliAdams}
H. D. Sherali and W. P. Adams, \newblock A reformulation-linearization technique for solving
discrete and continuous nonconvex problems. Kluwer, Dordrecht 1998.

\bibitem{Fraticelli}
H. D. Sherali and B. M. P. Fraticelli, \newblock Enhancing RLT
relaxations via a new class of semidefinite cuts. {\em J. Global
  Optim.} 22, pp.~233-261, 2002.

\bibitem{Shor}
N.Z. Shor, \newblock Quadratic optimization problems. {\em Tekhnicheskaya Kibernetika}, 1, 1987.

\bibitem{KrishnanMitchell1}
K. Sivaramakrishnan and J. Mitchell, \newblock Semi-infinite linear programming approaches to semidefinite programming (SDP) problems. Novel Approaches to Hard Discrete Optimization, edited by P.M. Pardalos and H. Wolkowicz, {\em Fields Institute Communications Series, American Math. Society}, 2003.

\bibitem{KrishnanMitchell2} 
\bysame, \newblock Properties of a cutting plane method for semidefinite programming, Technical Report, Department of Mathematics, North Carolina State University, September 2007.

\bibitem{SDPT3}
K. C. Toh, M. J. Todd and R. H. T\"{u}t\"{u}nc\"{u}, \newblock SDPT3:
A MATLAB software for semidefinite-quadratic-linear
programming. Available at\\
\href{http://www.math.nus.edu.sg/~mattohkc/sdpt3.html}
{\tt\small{http://www.math.nus.edu.sg/\~{}mattohkc/sdpt3.html}}

\bibitem{VandenbergheBoyd} 
L. Vandenberghe and S. Boyd, \newblock Semidefinite Programming. {\em
  SIAM Review} 38 (1), pp.~49-95, 1996.

\bibitem{Vandenbussche} 
D. Vandenbussche and G. L. Nemhauser, \newblock A branch-and-cut
algorithm for nonconvex quadratic programs with box constraints. {\em
  Math. Prog.} 102(3), pp.~559-575, 2005.

\bibitem{Wolkowicz}
H. Wolkowicz, R. Saigal and L. Vandenberghe, \newblock Handbook of Semidefinite Programming: Theory, Algorithms, and Applications. Springer, 2000.
\end{thebibliography}
\end{document}